# Comment

**Trevor Hastie and Ji Zhu**

We congratulate the authors for a well written and thoughtful survey of some of the literature in this area. They are mainly concerned with the geometry and the computational learning aspects of the support vector machine (SVM). We will therefore complement their review by discussing from the statistical function estimation perspective. In particular, we will elaborate on the following points:

- Kernel regularization is essentially a generalized ridge penalty in a certain feature space.
- In practice, the effective dimension of the data kernel matrix is not always equal to $n$, even when the implicit dimension of the feature space is infinite; hence, the training data are not always perfectly separable.
- Appropriate regularization plays an important role in the success of the SVM.
- The SVM is not fundamentally different from many statistical tools that our statisticians are familiar with, for example, penalized logistic regression.

We acknowledge that many of the comments are based on our earlier paper Hastie, Rosset, Tibshirani and Zhu (2004).

## KERNEL REGULARIZATION AND THE GENERALIZED RIDGE PENALTY

Given a positive definite kernel $K(\mathbf{x}, \mathbf{x}')$, where $\mathbf{x}, \mathbf{x}'$ belong to a certain domain $\mathcal{X}$, we consider the general function estimation problem

$$(1) \quad \min_{\beta_0, f} \sum_{i=1}^n \ell(y_i, \beta_0 + f(\mathbf{x}_i)) + \frac{\lambda}{2} \|f(\mathbf{x})\|_{\mathcal{H}_K}^2.$$

*Trevor Hastie is Professor, Department of Statistics, Stanford University, Stanford, California 94305, USA e-mail: hastie@stanford.edu. Ji Zhu is Assistant Professor, Department of Statistics, University of Michigan, Ann Arbor, Michigan 48109, USA e-mail: jizhu@umich.edu.*

Here $\ell(\cdot, \cdot)$ is a convex loss function that describes the "closeness" between the observed data and the fitted model, and $f$ is an element in the span of $\{K(\cdot, \mathbf{x}'), \mathbf{x}' \in \mathcal{X}\}$. More precisely, $f \in \mathcal{H}_K$ is a function in the reproducing kernel Hilbert space $\mathcal{H}_K$ (RKHS) generated by $K(\cdot, \cdot)$ (see Burges, 1998; Evgeniou, Pontil and Poggio, 2000; and Wahba, 1999, for details).

Suppose the positive definite kernel $K(\cdot, \cdot)$ has a (possibly finite) eigenexpansion,

$$K(\mathbf{x}, \mathbf{x}') = \sum_{j=1}^\infty \delta_j \phi_j(\mathbf{x}) \phi_j(\mathbf{x}'),$$

where $\delta_1 \geq \delta_2 \geq \cdots \geq 0$ are the eigenvalues and $\phi_j(\mathbf{x})$'s are the corresponding eigenfunctions. Elements of $\mathcal{H}_K$ have an expansion in terms of these eigenfunctions

$$(2) \quad f(\mathbf{x}) = \sum_{j=1}^\infty \beta_j \phi_j(\mathbf{x}),$$

with the constraint that

$$\|f\|_{\mathcal{H}_K}^2 \stackrel{\text{def}}{=} \sum_{j=1}^\infty \beta_j^2 / \delta_j < \infty,$$

where $\|f\|_{\mathcal{H}_K}$ is the norm induced by $K(\cdot, \cdot)$.

Then we can rewrite (1) as

$$(3) \quad \min_{\beta_0, \boldsymbol{\beta}} \sum_{i=1}^n \ell\left(y_i, \beta_0 + \sum_{j=1}^\infty \beta_j \phi_j(\mathbf{x}_i)\right) + \lambda \sum_{j=1}^\infty \frac{\beta_j^2}{\delta_j},$$

and we can see that the regularization term $\|f\|_{\mathcal{H}_K}^2$ in (1) can be interpreted as a generalized ridge penalty, where eigenfunctions with small eigenvalues in the expansion (2) get penalized more and vice versa.

Formulation (3) seems to be an infinite dimensional optimization problem, but according to the representer theorem (Kimeldorf and Wahba, 1971; Wahba 1990), the solution is finite dimensional and has the form

$$f(\mathbf{x}) = \sum_{i=1}^n \alpha_i K(\mathbf{x}, \mathbf{x}_i).$$

Using the reproducing property of $\mathcal{H}_K$, that is, $\langle K(\cdot, \mathbf{x}_i), K(\cdot, \mathbf{x}_{i'}) \rangle = K(\mathbf{x}_i, \mathbf{x}_{i'})$, (3) also reduces to a finite-dimensional criterion,

$$(4) \qquad \min_{\beta_0,\alpha} L(\mathbf{y}, \beta_0 + \mathbf{K}\alpha) + \lambda \alpha^T \mathbf{K}\alpha.$$

Here we use vector notation, $\mathbf{K}$ is the $n \times n$ data kernel matrix with elements equal to $K(\mathbf{x}_i, \mathbf{x}_{i'}), i, i' = 1, \ldots, n$, and $L(\mathbf{y}, \beta_0 + \mathbf{K}\alpha) = \sum_{i=1}^n \ell(y_i, \beta_0 + f(\mathbf{x}_i))$. We reparametrize (4) using the eigendecomposition of $\mathbf{K} = \mathbf{U}\mathbf{D}\mathbf{U}^T$, where $\mathbf{D}$ is diagonal and $\mathbf{U}$ is orthogonal. Let $\mathbf{K}\alpha = \mathbf{U}\beta$, where $\beta = \mathbf{D}\mathbf{U}^T\alpha$. Then (4) becomes

$$(5) \qquad \min_{\beta_0,\beta} L(\mathbf{y}, \beta_0 + \mathbf{U}\beta) + \lambda \beta^T \mathbf{D}^{-1} \beta.$$

Now the columns of $\mathbf{U}$ are unit-norm basis functions that span the column space of $\mathbf{K}$; and again, we see that those members that correspond to small eigenvalues (the elements of the diagonal matrix $\mathbf{D}$) get heavily penalized and vice versa.

In the machine learning community, people tend to view the kernel as providing an implicit map of $\mathbf{x}$ from $\mathcal{X}$ to a certain high-dimensional feature space, and $K(\cdot, \cdot)$ computes inner products in this (possibly infinite-dimensional) feature space. Specifically, the features are

$$h_j(\mathbf{x}) = \sqrt{\delta_j}\phi_j(\mathbf{x}) \quad \text{or} \quad \mathbf{h}(\mathbf{x}) = (h_1(\mathbf{x}), h_2(\mathbf{x}), \ldots)^T,$$

and we have

$$K(\mathbf{x}, \mathbf{x}') = \langle \mathbf{h}(\mathbf{x}), \mathbf{h}(\mathbf{x}') \rangle.$$

Furthermore, let $\theta_j = \beta_j/\sqrt{\delta_j}$ and $\mathbf{H} = \mathbf{U}\mathbf{D}^{1/2}$. Then (3) and (5) become

$$(6) \qquad \min_{\beta_0,\boldsymbol{\theta}} \sum_{i=1}^n \ell\left(y_i, \sum_{j=1}^\infty \theta_j h_j(\mathbf{x})\right) + \lambda \sum_{j=1}^\infty \theta_j^2$$

and

$$(7) \qquad \min_{\beta_0,\boldsymbol{\theta}} L(\mathbf{y}, \beta_0 + \mathbf{H}\boldsymbol{\theta}) + \lambda \boldsymbol{\theta}^T \boldsymbol{\theta},$$

respectively. This shows kernel regularization as an exact ridge penalty in the feature space, but unlike (3) and (5), it hides the fact that eigenfunctions are differentially penalized according to their corresponding eigenvalues.

To illustrate the point, we consider a simple example. The data $x_i$'s are one-dimensional and were generated from the standard Gaussian distribution ($n = 50$). The radial kernel function $K(x, x') = \exp(-\gamma \|x - x'\|^2)$ was used, with $\gamma = 1$. Figure 1 shows the eigenvalues of the kernel matrix $\mathbf{K}$. The left panel of Figure 2 shows the first 16 eigenvectors of $\mathbf{K}$ (columns of $\mathbf{U}$) and the right panel shows the corresponding features (columns of $\mathbf{H}$). As we can see from the left panel of Figure 2, eigenvectors with large eigenvalues (hence get penalized less) tend to be smooth, while eigenvectors with small eigenvalues (hence get penalized more) tend to be wiggly; therefore, $\|f\|^2_{\mathcal{H}_K}$ is also always interpreted as a roughness measure of the function $f$. We can also see from the right panel of Figure 2 that many of the features are "norm challenged," that is, they are squashed down dramatically by their eigenvalues.

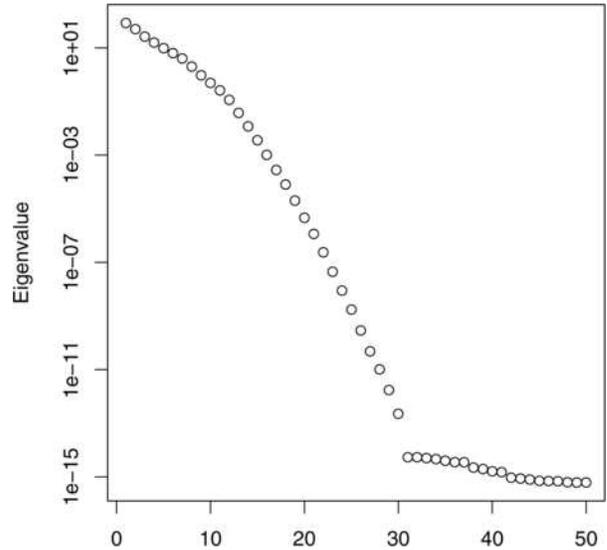

FIG. 1. *Eigenvalues (on the log scale) of the data kernel matrix $\mathbf{K}$.*

## EFFECTIVE DIMENSION OF THE DATA KERNEL MATRIX

As we have seen in the previous section, the kernel $K(\cdot, \cdot)$ maps $\mathbf{x}$ from its original input space to some high-dimensional feature space $\mathbf{h}(\mathbf{x})$. In the case of classification, it is sometimes argued that the implicit feature space can be infinite-dimensional (e.g., via the radial basis kernel), which suggests that perfect separation of the training data is always possible. However, this is not always true in practice.

To illustrate the point, we consider a two-class classification example. The data were generated from a pair of mixture Gaussian densities, described in detail by Hastie, Tibshirani and Friedman (2001). The radial kernel function $K(\mathbf{x}, \mathbf{x}') = \exp(-\gamma \|\mathbf{x} - \mathbf{x}'\|^2)$ was used. Four different values of $\gamma$ (0.1, 0.5, 1 and 5) were tried. For each of the values of $\gamma$, the SVM was fitted for a sequence of values of $\lambda$, ranging from the most regularized model to the least regularized model.

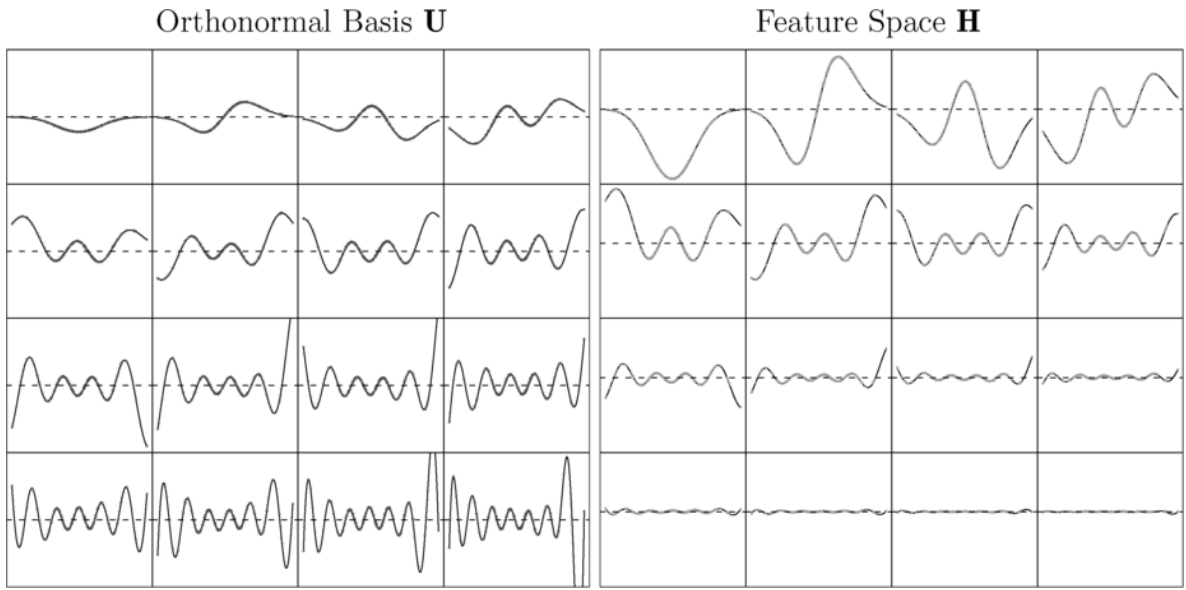

FIG. 2. *The left panel shows the eigenvectors of the data kernel matrix* **K** *and the right panel shows the corresponding features (eigenvectors scaled by the corresponding eigenvalues).*

TABLE 1
*Results for the mixture simulation example*

| $\gamma$ | 5 | 1 | 0.5 | 0.1 |
|---|---|---|---|---|
| Training errors | 0 | 12 | 21 | 33 |
| Effective rank | 200 | 177 | 143 | 76 |

We were at first surprised to discover that not all these sequences achieved zero training errors on the 200 training data points, at their least regularized fit. The minimal training errors and the corresponding values for $\gamma$ are summarized in Table 1. The second row of the table shows the effective rank of the data kernel matrix **K** (which we defined to be the number of eigenvalues greater than $10^{-12}$). This $200 \times 200$ matrix has elements $K(\mathbf{x}_i, \mathbf{x}_{i'}), i, i' = 1, \ldots, 200$. In this example, a full rank **K** is required to achieve perfect separation. Similar observations have also appeared in Williams and Seeger (2000) and Bach and Jordan (2002).

This emphasizes again the fact that not all features in the feature map implied by $K(\cdot, \cdot)$ are of equal stature (see the right panel in Figure 2); many of them are shrunken way down to zero. The regularization in (3) and (5) penalizes unit-norm eigenvectors by the inverse of their eigenvalues, which effectively annihilates some, depending on $\gamma$. Small $\gamma$ implies wide, flat kernels and a suppression of wiggly, rough functions. Figure 3 shows the eigenvalues of **K** for the four values of $\gamma$. The larger eigenvalues correspond in this case to smoother eigenfunctions, the smaller ones to rougher. The rougher eigenfunctions get penalized exponentially more than the smoother ones. Hence, for smaller values of $\gamma$, the effective dimension of the feature space is truncated.

## THE NEED FOR CAREFUL REGULARIZATION

The SVM has been very successful for the classification problem and gained a lot of attention in the machine learning community in the past ten years. Many papers have been published to explain why it performs so well. Most of this literature concentrates on the concept of margin. Various misclassification error bounds have been derived based on the margin (Vapnik, 1995; Bartlett and Shawe-Taylor, 1999; Shawe-Taylor and Cristianini, 1999).

However, our view is a little different from that based on the concept of margin. Several researchers have noted the relationship between the SVM and regularized function estimation in RKHS (Evgeniou, Pontil and Poggio, 2000; Wahba, 1999). The regularized function estimation problem contains two parts: a loss function and a penalty term [e.g., (1)]. The SVM uses the so-called hinge loss function (see Figure 5). The margin maximizing property of the SVM derives from the hinge loss function. Hence margin maximization is by nature a nonregularized objective, and solving it in high-dimensional space is likely to lead to overfitting and bad prediction performance. This has been observed in practice by many researchers, in particular Breiman (1999) and Marron and Todd (2002).

The *loss + penalty* formulation emphasizes the role of regularization. In many situations we have sufficient features (e.g., gene expression arrays) to guarantee separation. We may nevertheless avoid the maximum margin separator ($\lambda \downarrow 0$) in favor of a more regularized solution.

Figure 4 shows the test error as a function of $\lambda$ for the mixture data example. Here we see a dramatic range in the correct choice of $\lambda$. When $\gamma = 5$, the most regularized model is called for. On the other hand, when $\gamma = 0.1$, we would want to choose among the least regularized models. Depending on the value of $\gamma$, the optimal $\lambda$ can occur at either end of the spectrum or anywhere in between, emphasizing the need for careful selection.

## CONNECTION WITH OTHER STATISTICAL TOOLS

Last, we would like to comment on the connection between the SVM and some statistical tools that statisticians are familiar with.

As we have seen in previous sections, what is special with the SVM is not the regularization term, but is rather the loss function, that is, the hinge loss. Lin (2002) pointed out that the hinge loss is Bayes consistent, that is, the population minimizer of the loss function agrees with the Bayes rule in terms of classification. This is important in explaining the success of the SVM, because it implies that the SVM is trying to implement the Bayes rule.

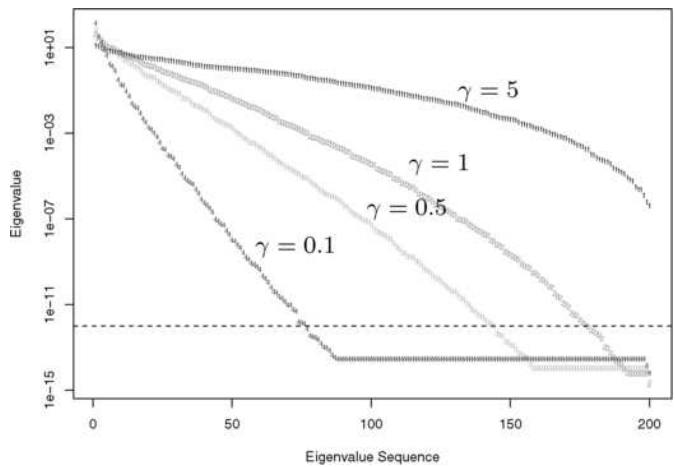

FIG. 3. *The eigenvalues (on the log scale) for the data kernel matrices* **K** *that correspond to the four values of* $\gamma$.

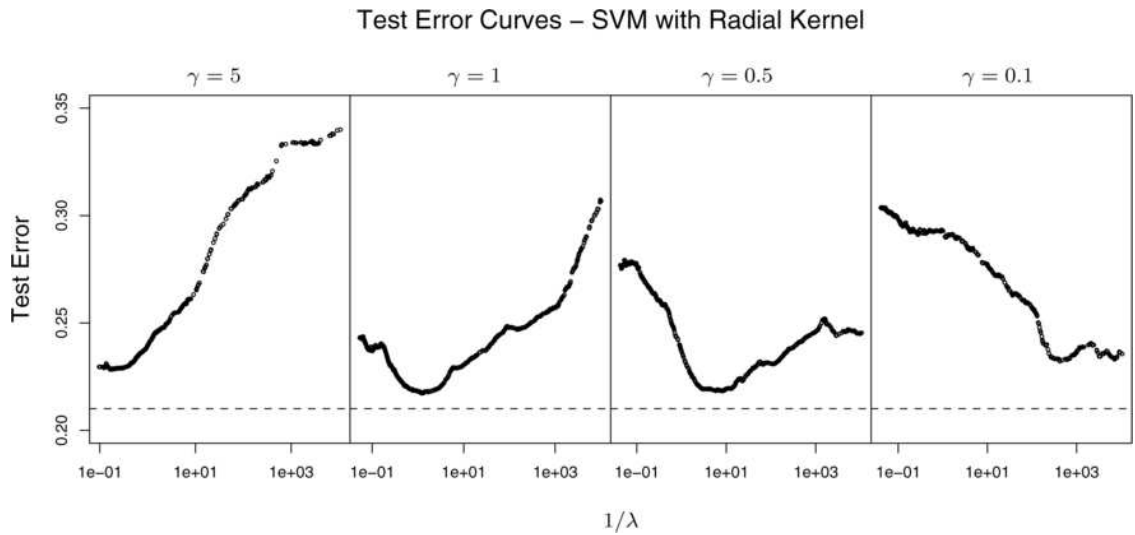

FIG. 4. *Test error curves for the mixture example, using four different values for the radial kernel parameter* $\gamma$. *Large values of* $\lambda$ *correspond to heavy regularization, small values of* $\lambda$ *to light regularization. Depending on the value of* $\gamma$, *the optimal* $\lambda$ *can occur at either end of the spectrum or anywhere in between, emphasizing the need for careful selection.*

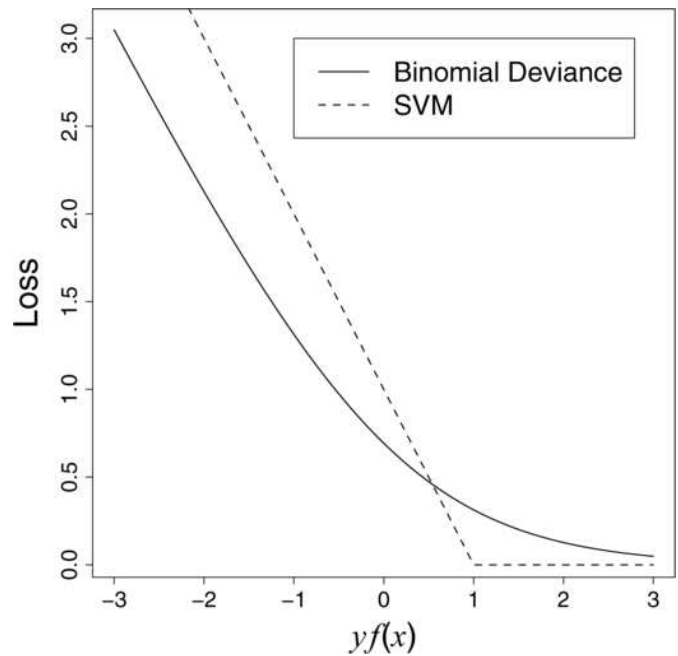

Fig. 5. *Comparing the hinge loss and the binomial deviance, $y \in \{-1, 1\}$.*

On the other hand, notice that the hinge loss and the binomial deviance have very similar shapes (see Figure 5): both increase linearly as $yf$ gets very small (negative), and both encourage $y$ and $f$ to have the same sign. Hence it is reasonable to conjecture that by replacing the hinge loss of the SVM with the binomial deviance, which is also Bayes consistent, we should be able to get a fitted model that performs similarly to the SVM. In fact, in Zhu and Hastie (2005), we show that under certain conditions, the classification boundary of the resulting penalized logistic regression (using the binomial deviance) and that of the SVM coincide. Penalized logistic regression has been studied by many statisticians (see Green and Silverman, 1994; Wahba, Gu, Wang and Chappell, 1995; and Lin et al., 2000, for details). We understand why it can work well. The same reasoning could be applied to the SVM.

Penalized logistic regression is not the only model that performs similarly to the SVM; replacing the hinge loss with any sensible loss function will give a similar result, for example, the exponential loss function of boosting (Freund and Schapire, 1997) and the squared error loss (Zhang and Oles, 2001; Bühlmann and Yu, 2003). These loss functions are all Bayes consistent. The binomial deviance and the exponential loss are margin-maximizing loss functions, but the squared error loss is not. The distance weighted discrimination (Marron and Todd, 2002) is designed specifically for *not* maximizing the margin and works well with high-dimensional data, which in a way also justifies that margin maximization is not the key to the success of the SVM.

## ACKNOWLEDGMENTS

Hastie is partially supported by NSF Grant DMS-05-05676. Zhu is partially supported by NSF Grant DMS-05-05432.